# Stabilization of linear time-delayed systems by delayed feedback

Mohammad Mousa-Abadian, Sayed Hodjatollah Momeni-Masuleh and Mohammad Haeri

Advanced Control Systems Lab, Electrical Engineering, Sharif University of Technology, Tehran, Iran

Email: m.abadian@shahed.ac.ir, momeni@shahed.ac.ir, haeri@sina.sharif.edu

**Abstract:** In this paper, stability analysis of time delay systems is considered based on decomposition of the systems to subsystems. The decomposition process needs matrices of these systems to be simultaneously block triangularize. We show that a finite set of matrices are simultaneously block triangularize if and only if they have a common invariant subspace. The importance of the decomposition is highlighted by some systems that their characteristic equation has repeated roots on the imaginary axis. For such systems, the cluster treatment method and the direct method cannot be employed to analyze the stability of them, unless, we decompose the systems to subsystems. On the other hand, it has been shown that stabilization of time delay systems can be done by delayed feedback. Moreover, this kind of feedback is applied to design a controller.

## 1. Introduction

System of linear delay differential equations (DDEs) appears naturally in many branches of science and engineering. A common way of describing these equations is the state-space representation which widely used in control theory. Unlike ordinary differential equations (ODEs), in DDEs the rate of change of a time-dependent process not only depends on the current state, but also depends on its history state.

The stability analysis of such systems is vital from theoretical and practical point of view [1-3]. Due to the presence of the delayed term, the characteristic equation is a quasi-polynomial instead of a polynomial. Therefore, a complicated stability analysis is necessitated. Rekasius [4], Walton and Marshall [5] and Olgac and Sipahi [6] are addressed some methods for the stability analysis of DDEs. However, these methods cannot be applied for an arbitrary system of DDEs. Mesbahi and Haeri [7] presented an example in which the characteristic equation has a multiple root with multiplicity greater than one and showed that the method of Olgac and Sipahi [6] cannot be used to analyze the stability. They removed repeated roots by decomposing original system to several subsystems with lower dimension. One of the most common types of DDEs is retarded DDE. In a retarded linear system of DDEs the derivative term $\dot{x}$, doesnot depend on delay. In this paper, we present some examples that the cluster treatment [6] for retarded DDEs cannot be applied for analyzing the stability of time delay systems and the direct method [5] leads to an ambiguous result in detecting the number of unstable poles. Meanwhile, decomposing the original system to subsystems resolves these difficulties.

To decompose a single linear delay system that has two matrices, one needs an invertible transformation such that simultaneously transforms both matrices to a block triangular or diagonal form.

Recently, because of applications of simultaneous triangularization and simultaneous diagonalization in multidimensional systems [8], discrete time switching systems [9] and differential inclusions [10-11], researches on this topic has been widely increased.

Two questions are crucial in simultaneous triangularization (or in simultaneous diagonalization): First, when two or more generally, a finite set of matrices can be transformed simultaneously to a block triangular (or diagonal) form? Second, what kind of transformation can put matrices into a block triangular (or diagonal) form simultaneously? It is easy to show that every set of commutative matrices can be simultaneously transformed in an upper triangular form, but the converse is not true [12-13]. One of the most famous classical theorems is McCoy's theorem [14] that provides necessary and sufficient conditions in response to the first question. Levitzky [13] proved that a semigroup of nilpotent matrices can be put simultaneously into a triangular form. Radjavi [15] provided a trace condition which is equivalent to simultaneous triangularization. Uhlig [16] presented a necessary and sufficient condition to converting two real symmetric matrices into simultaneous block diagonalization form. Shapiro [17] gives a theorem to put a set of square matrices into the simultaneous block upper triangular form. In Ref. [19], the problem of simultaneous triangularization of matrices is considered only for low rank cases and the nonderogatory case. Dubi [20] proposed an algorithm to construct a simultaneous triangularization of a set of matrices. He presented an algorithm that answers the first and second question for non-block triangularization. His algorithm uses Shemesh's idea [21] to compute the invertible transformation. Kaczorek [22] proved a necessary and sufficient condition for block simultaneous triangularization of a set of matrices. To construct the invertible transformation, it requires a full column rank matrix that cannot be computed by an algorithmic way.

However, the presence of time delay in a system can cause various difficulties, e.g. it can be a source for instability in the system; nevertheless it can be useful in some cases. Pyragas [23] applied delayed feedback to control chaos and also he employed this feedback to stabilize the unstable periodic orbits.

Usually, one of the major goals in control theory is controlling an equilibrium solution or equivalently, the regulator problem. In fact, in a regulator problem, one needs to obtain an asymptotically stable steady state solution which attracts all nearby initial conditions. Dahms et al. [24] considered the extended time-delayed feedback method to control unstable steady states.

Usually time delay is a source for instability. Nevertheless, Abdallah et al. [25] showed that some oscillatory systems can be stabilized by a delayed feedback. On the other hand, it has been shown that the performance of a system can be improved by a delayed feedback and also disturbance attenuation and robustness against parameters variation can be modified [26]. Moreover, Kwon et al. [26] obtained the state feedback tracking controller by the delayed feedback method.

In this paper, we show that the block simultaneous triangularization of a set of square matrices is equivalent to existence of a common invariant subspace for them. We provide a lemma that characterizes the invariant subspaces of a matrix by means generalized eigenvectors [27]. Furthermore, we present some linear time delay systems that their stability analysis cannot be investigated by the direct method [5] and the cluster treatment method [6]. Also, we show that an unstable time-delay system can be stabilized by the delayed feedback method. In addition, we adopt this feedback to put poles of a time-delay system in certain places of complex plane and so a desired response can be obtained. More precisely, by the delayed feedback method, the settling time of the system is remarkably reduced.

The remaining of the paper is organized as follows. In Section 2, we introduce some required mathematical details and problem statement. Stability analysis of the linear time delay systems using decomposition of matrices is considered in Section 3. In

Section 4 stabilization of an unstable time-delay system by delayed feedback method is discussed. Section 5 is devoted to design a controller via the delayed feedback. Finally, conclusion is available in Section 6.

## 2. Mathematical details and problem statement

2. 1 Definitions, lemmas and theorems

Let us to introduce some notations that are used throughout this paper. $\mathbb{R}^n$ and $\mathbb{C}^n$ are real and complex Euclidean spaces respectively. $\mathbb{R}^{n \times n}$ ($\mathbb{C}^{n \times n}$) is space of real (complex) square matrices. Also, we denote crossing frequencies and their corresponding delays by $\omega_{ck}$ and $\tau_{kl}$ respectively, where $k = 1,2,\ldots,n$, $l = 1,2,\ldots$. $F(s, \tau)$ denotes characteristic equation of a time-delay system.

In this part of the article, first, we explain some definitions, theorems and lemmas that are useful for simultaneous block triangularization of a set of square matrices. Then, the controllability theorem for the time delay systems is expressed.

**Definition 1.** Let $A_1, A_2, \ldots, A_N$ be a set of square matrices in $\mathbb{R}^{n \times n}$, where $N$ is a positive integer number. These matrices are said to be simultaneously block triangularize with dimension $k$, if there exists an invertible transformation $T$ such that

$$TA_i T^{-1} = \tilde{A}_i = \begin{bmatrix} \tilde{A}_{i1} & \tilde{A}_{i2} \\ 0 & \tilde{A}_{i4} \end{bmatrix}, \quad i = 1, 2, \ldots, N, \tag{1}$$

where $\tilde{A}_{i1} \in \mathbb{R}^{k \times k}, \tilde{A}_{i2} \in \mathbb{R}^{k \times (n-k)}, \tilde{A}_{i4} \in \mathbb{R}^{(n-k) \times (n-k)}, i = 1, \ldots, N$, and $1 \leq k < n$.

**Example 1 [22].** Consider the following matrices

$$A_1 = \begin{bmatrix} 1 & 1 & 0 \\ 0 & 2 & 0 \\ 0 & 3 & 1 \end{bmatrix}, A_2 = \begin{bmatrix} 0 & 1 & 0 \\ 0 & 4 & 0 \\ 2 & 2 & 0 \end{bmatrix}.$$

If we put

$$T = \begin{bmatrix} 1 & 0 & 0 \\ 0 & 0 & 1 \\ 0 & 1 & 0 \end{bmatrix},$$

then we get the following

$$TA_1T^{-1} = \begin{bmatrix} 1 & 0 & 1 \\ 0 & 1 & 3 \\ 0 & 0 & 2 \end{bmatrix}, TA_2T^{-1} = \begin{bmatrix} 0 & 3 & 1 \\ 2 & 0 & 2 \\ 0 & 0 & 4 \end{bmatrix}.$$

So, in this example $k = 2$.

**Definition 2.** A vector subspace $V \subset \mathbb{R}^n$ is said to be $(A_1, A_2, \ldots, A_N)$-invariant if $A_i v \in V \ \forall v \in V, i = 1, 2, \ldots, N$.

**Definition 3 [12].** Let $A$ be a square matrix that belongs to $\mathbb{R}^{n \times n}$. $x_0, x_1, \ldots, x_k$ is called a Jordan chain of $A$ corresponding to $\lambda_0$ if $x_0 \neq 0$ and the following relations hold

$$Ax_0 = \lambda_0 x_0,$$
$$Ax_1 - \lambda_0 x_1 = x_0,$$
$$Ax_2 - \lambda_0 x_2 = x_1,$$
$$\vdots$$
$$Ax_k - \lambda_0 x_k = x_{k-1}.$$

The first equation (together with $x_0 \neq 0$) means that $x_0$ is an eigenvector of $A$ corresponding to $\lambda_0$. The vectors $x_0, x_1, \ldots, x_k$ are called generalized eigenvector of $A$ corresponding to the eigenvalue $\lambda_0$ and the eigenvector $x_0$.

**Definition 4[27].** Let $A \in \mathbb{R}^{n \times n}$ and define a set of Jordan vectors for $A$ to be a set of linearly independent vectors in $\mathbb{C}^n$ made up of a union of Jordan chains.

In the following theorem we show that simultaneous block triangularization of a finite set of matrices is equivalent to existence of a common invariant $k$ dimensional subspace for them. For the sake of convenience, we express this theorem only for the set of two matrices.

**Definition 5[1].** Consider the following system

$$\dot{x}(t) = \sum_{k=0}^{N} A_k x(t - \tau_k) + Bu(t),$$

where $x(t) \in \mathbb{R}^n$, $u(t) \in \mathbb{R}^n$ and $\tau_0 = 0$. The system is $\mathbb{R}^n$- controllable on $[t_0, t_1]$ if $\forall x_0 \in C[-\tau_N, 0]$ and $\forall x_1 \in \mathbb{R}^n$ there exists piecewise-continuous $u(t) = u(t, x_0, x_1)$ such that the solution of the system with the initial condition $x_{t_0} = x_0$ satisfies $x_{t_1} = x_1$.

**Theorem 1.** The matrices $A_1, A_2 \in \mathbb{R}^{n \times n}$ are simultaneously block triangularize with dimension $k$, if and only if there exists an $(A_1, A_2)$-invariant $k$-dimensional subspace $W \subset \mathbb{R}^n$.

**Proof.** Let $W \subset \mathbb{R}^n$ be an arbitrary $k$-dimensional vector subspace and $Q = [Q_k \ Q_{n-k}] = [q_1 \ \cdots \ q_k \ q_{k+1} \ \cdots \ q_n]$ be a nonsingular matrix where the first k columns, i.e. $Q_k = [q_1 \ \cdots \ q_k]$, of $Q$ form a basis for the subspace $W$. By assuming

$$Q^{-1} A_i Q = \begin{bmatrix} \tilde{A}_{i1} & \tilde{A}_{i2} \\ \tilde{A}_{i3} & \tilde{A}_{i4} \end{bmatrix}, \quad i = 1,2, \tag{2}$$

we show that $\tilde{A}_{13} = \tilde{A}_{23} = 0$ iff the first $k$ columns of $Q$ form a basis for the subspace $W$. To do this end, let $[a_1^{i,1} \ \cdots \ a_k^{i,1}]$ be the $k$ columns of $\tilde{A}_{i1}$ and $[a_1^{i,3} \ \cdots \ a_k^{i,3}]$ be the $k$ columns of $\tilde{A}_{i3}$ for $i = 1,2$, respectively. Since $Q$ is a nonsingular matrix, by the Definition 1, we have

$$A_i [Q_k \ Q_{n-k}] = [Q_k \ Q_{n-k}] \begin{bmatrix} \tilde{A}_{i1} & \tilde{A}_{i2} \\ \tilde{A}_{i3} & \tilde{A}_{i4} \end{bmatrix}, \quad i = 1,2,$$

$$= [Q_k \tilde{A}_{i1} + Q_{n-k} \tilde{A}_{i3} \ Q_k \tilde{A}_{i2} + Q_{n-k} \tilde{A}_{i4}], \tag{3}$$

Therefore, by equating corresponding columns in (3), we obtain the following relations

$$A_i q_j = [Q_k \tilde{A}_{i1} + Q_{n-k} \tilde{A}_{i3}]_j, \quad j = 1, \ldots, k, i = 1,2.$$

So, for $j = 1, \ldots, k, i = 1,2$, we have

$$A_i q_j = Q_k a_j^{i,1} + Q_{n-k} a_j^{i,3}$$

As the first $k$ columns of $Q$ are linearly independent, therefore, $\tilde{A}_{13} = \tilde{A}_{23} = 0$ which means that for $i = 1,2$, $W$ is $A_i$-invariant. Clearly, if $W$ is $A_i$-invariant, then $\tilde{A}_{13} = \tilde{A}_{23} = 0$. This completes the proof. □

**Remark 1**. It is clear that $Q$ is not unique.

Following corollaries are immediate consequence of Theorem 1.

**Corollary 1.** If $A_1, A_2, \ldots, A_N$ have $k$ common eigenvectors, then they can be transformed simultaneously into a $k$-dimensional block triangular form.

**Corollary 2.** Let $n = 2$. $A_1, A_2, \ldots, A_N \in \mathbb{R}^{n \times n}$ are simultaneously block triangular form if and only if they have a common eigenvector.

Now, Theorem 1 leads us to the following theorem [22]. In other words, the following theorem is a consequence of Theorem 1.

**Theorem 2** [22]. Let $A_1, A_2, \ldots, A_N$ be square matrices in $\mathbb{R}^{n \times n}$, These matrices can be put simultaneously in form (1) by means of transformation $T$, if and only if there exists a full column rank matrix $J \in \mathbb{R}^{n \times r}$ such that rank $[J \quad A_i J] = r$ for $i = 1, 2, \ldots, N$.

The following proposition characterizes $k$ dimensional invariant subspaces for a square matrix $A$.

**Proposition 1**. Let $M$ be a real $n$ by $n$ matrix. A $k$ dimensional subspace $W \subset \mathbb{R}^n$ is $A$-invariant iff $W$ has a basis consisting of a set of Jordan vectors for $M$.

**Proof.** Assume $W$ has a set of Jordan vectors, say $\{x_1, x_2, \ldots, x_k\}$, for $M$ as a basis. By the assumption, $W = \text{span} < x_1, \ldots, x_k >$. Since $\{x_1, x_2, \ldots, x_k\}$ belongs to a

Jordan chain, so by proposition 1.3.1 in Ref. [12], $W$ is $A$-invariant.

Conversely, let $X = [x_1 \quad x_2 \quad \ldots \quad x_k]$ be a $n \times k$ matrix whose columns form an arbitrary basis for W. Since W is A-invariant, there exists $G \in \mathbb{R}^{k \times k}$ such that $MX = XG$. The Jordan matrix decomposition of $G$ can be written as $G = SJS^{-1}$ for some $S$, which leads us to $MXS = XSJ$ and therefore, we get $J = (XS)^{-1}M(XS)$. Here, in fact, the columns of the matrix $XS$ form the Jordan vector for $M$. □

Now, we give the following theorem that is crucial for controllability of the time delay systems.

**Theorem 3 [28]**. If $(A_0 + A_1, B)$ is controllable, then the following system is controllable

$$\dot{x}(t) = A_0 x(t) + A_1 x(t-\tau) + Bu(t), \quad A_0, A_1, B \in \mathbb{R}^{n \times n}, u \in \mathbb{R}^{n \times 1}.$$

**2. 2 Problem statement**

In this paper, first, we consider the stability analysis of the following linear time-delay system:

$$\dot{x}(t) = A_1 x(t) + A_2 x(t-\tau), \quad A_1, A_2 \in \mathbb{R}^{n \times n},$$

where $\tau > 0$ is time delay. We will show that stability analysis of some systems required to decompose them to subsystems with lower dimension and then analyze the stability of each subsystem and finally stability of the whole system is achieved.

Then, by the delayed feedback method we attempt to stabilize the following system

$$\dot{x}(t) = A_0 x(t) + A_1 x(t-\tau) + Bu(t),$$

where $\tau > 0$ is a fixed delay and $A_0, A_1, B, u$ are the same as defined in theorem 3.

**3. Stability analysis of linear time delay systems via decomposition**

In this section we provide two examples that show the proposed method in [5] cannot recognize the number of unstable poles. In addition, the cluster treatment [6] also fails to analyze the stability of the systems.

**Example 2:** Consider the following system

$$\dot{x}(t) = A_1 x(t) + A_2 x(t-\tau), \tag{4}$$

where

$$A_1 = \begin{bmatrix} 3.2423 & -1.4176 & -2.7298 & 4.6267 \\ -1.0366 & -0.9812 & -0.7598 & -3.2319 \\ 2.0250 & 0.8723 & 0.0129 & 4.0908 \\ -0.9802 & 1.5668 & 1.2885 & -1.2741 \end{bmatrix},$$

$$A_2 = \begin{bmatrix} 1.4104 & 1.1252 & -0.1052 & 0.9652 \\ -0.2045 & -0.5965 & -0.2415 & -0.2683 \\ 0.4985 & 0.7644 & 0.1801 & 0.4498 \\ -0.3069 & 0.4843 & 0.4550 & 0.0060 \end{bmatrix}.$$

The characteristic equation of the system cross the imaginary axis at $s = \pm j$, $s = \pm\sqrt{3}j$ and corresponding delays are $\tau_k = \pi + 2\pi k$, $\tau_k = \frac{2k\pi}{\sqrt{3}}$, $k = 0,1,...$, respectively. Since all of roots of the characteristic equation are in the right half plane for $\tau = 0$, so the system is unstable without delay. By [5], since $\operatorname{sgn} W'(\omega^2)$ is positive at $\omega = \sqrt{3}$ and zero at $\omega = 1$, so the system is unstable for all $\tau$. In other words the method says all roots of the time-delay system (4) are in the right half plane. If we want to apply the method in [6], then we have $\frac{\partial F(s,\tau)}{\partial s} = 0, \frac{\partial F(s,\tau)}{\partial \tau} = 0$ at crossing frequency $\omega = 1$, and so the root tendency cannot be determined.

Now, we decompose the system as follows.

$A_1$ and $A_2$ have a common invariant subspace with dimension 2. A basis for this subspace is as following

$$E = \begin{bmatrix} 0.3878 & 0.8143 \\ -0.2562 & -0.1180 \\ 0.5371 & 0.2878 \\ -0.2094 & -0.1772 \end{bmatrix}.$$

In fact, linear combination of the columns of $E$ form a two dimensional $(A_1, A_2)$-

invariant subspace. One choice for transformation $T$ in (1) is

$$T = \begin{bmatrix} -1.0000 & 3.6667 & 4.3333 & 0 \\ 1.7321 & 1.6000 & 0 & 1.2500 \\ 0.2857 & 0 & 0.8333 & 2.6667 \\ 0 & 2.2500 & 1.3333 & 0.6667 \end{bmatrix}.$$

By applying this transformation to system (4), we get two subsystems which are as follow

$$\dot{z}_1(t) = \overline{A}_{11}z_1(t) + \overline{B}_{11}z_1(t-\tau), \quad (5)$$

$$\dot{z}_2(t) = \overline{A}_{22}z_2(t) + \overline{B}_{22}z_2(t-\tau), \quad (6)$$

where $x(t) = Tz(t)$ and

$$\overline{A}_{11} = \begin{bmatrix} 0 & 1 \\ -1 & 1 \end{bmatrix}, \overline{B}_{11} = \begin{bmatrix} 0 & 0 \\ 0 & 1 \end{bmatrix}, \overline{A}_{22} = \begin{bmatrix} 0 & 2 \\ -1 & 0 \end{bmatrix}, \overline{B}_{22} = \begin{bmatrix} 0 & 1 \\ 0 & 0 \end{bmatrix}.$$

Now, we employ the proposed method in [3] to each of system (5) and (6) separately. System (5) has a crossing frequency at $s = \pm j$ and $\text{sgn}\, W'(\omega^2)$ is zero for this frequency. So the system is unstable for all $\tau$. System (6) crosses the imaginary axis at $s = \pm\sqrt{3}j$, $s = \pm j$ and $\text{sgn}\, W'(\omega^2)$ at these frequency is positive and negative, respectively. Therefore, this system is stable for $\pi < \tau < \frac{2\pi}{\sqrt{3}}$. Thus, the system (4) has two stable poles for $\pi < \tau < \frac{2\pi}{\sqrt{3}}$. Roots of the characteristic equation of the system (5) and root locus for the system (6) are plotted in Fig. 1 and Fig. 2 respectively.

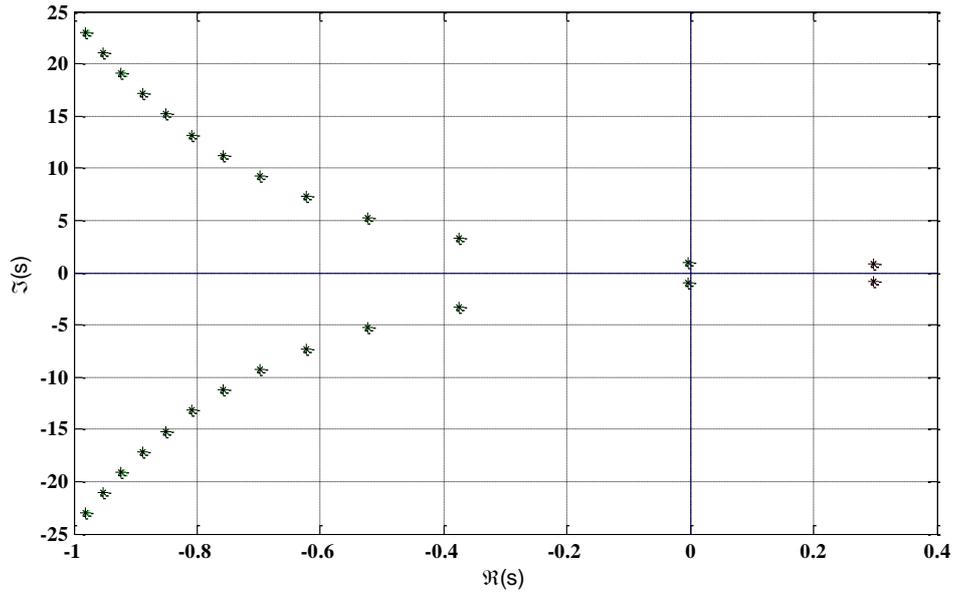

Fig. 1. Roots of the characteristic equation of the system (5), for $\tau = 3.2$.

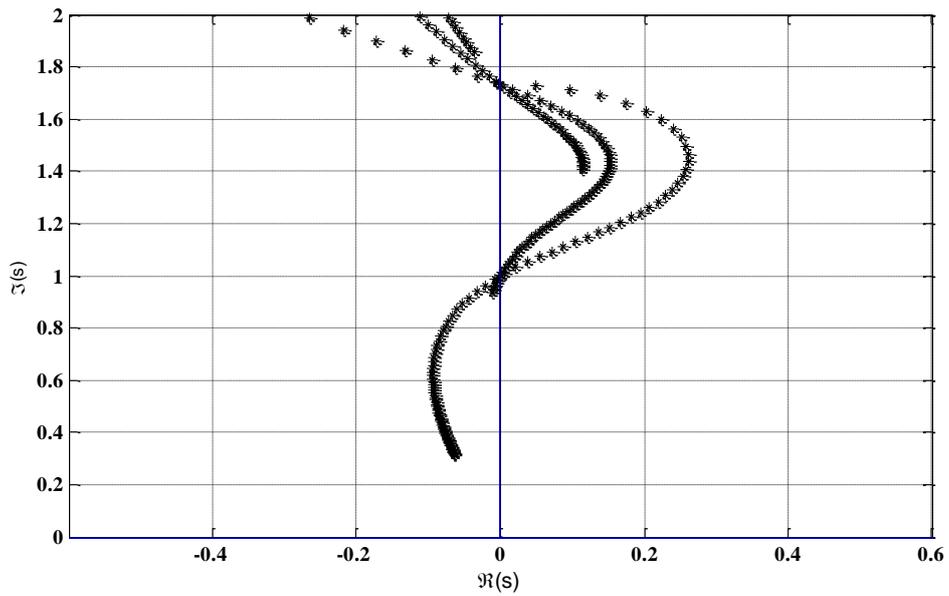

Fig. 2. The root locus of the system (6) near the imaginary axis, for $0 < \tau < 4$.

**Example 3:** Let us to analyze the stability of the following system

$$\dot{x}(t) = A_1 x(t) + A_2 x(t - \tau), \tag{7}$$

where

$$A_1 = \begin{bmatrix} -14.6102 & -4.9441 & 11.3503 & -11.5177 & -11.9699 \\ -3.9437 & -1.0804 & 3.4948 & -3.3674 & -3.2193 \\ 6.4695 & 0.5153 & -4.1521 & 3.9784 & 5.0394 \\ 6.0633 & 2.1406 & -4.6372 & 5.0694 & 4.8474 \\ 20.3590 & 4.5468 & -15.5102 & 13.5751 & 16.7733 \end{bmatrix},$$

$$A_2 = \begin{bmatrix} -11.1098 & -3.6577 & -2.2712 & -13.4823 & -4.0327 \\ -3.1263 & -1.0354 & -0.6680 & -3.7568 & -1.1390 \\ 4.8695 & 1.7361 & 1.6197 & 5.1076 & 1.8581 \\ 4.4403 & 1.4397 & 0.8037 & 5.5222 & 1.5967 \\ 16.3449 & 5.4846 & 3.8268 & 19.2118 & 6.0034 \end{bmatrix}.$$

First, we decompose the system and then we apply the method [5] to each subsystem.

The columns of the full column rank matrix

$$E = \begin{bmatrix} 1.2775 & -1.3977 \\ 0.6036 & -0.4111 \\ -0.5536 & 0.9967 \\ -0.5480 & 0.4946 \\ -1.9230 & 2.3550 \end{bmatrix},$$

Constitute a two dimensional $(A_1, A_2)$- invariant subspace.

We may choose the transformation $T$ in (1) as

$$T = \begin{bmatrix} 0.125 & 4.5 & 0.6667 & 2.75 & 0 \\ 1.4142 & 0.75 & 1.625 & 0 & 0.7071 \\ 3 & 0.8571 & 0 & 4.4286 & 1 \\ -1 & 0 & 2 & 2.6667 & -2 \\ 0 & 1.7321 & -1 & 1.4142 & 0.4286 \end{bmatrix}.$$

After applying this to system (7), two subsystems are available as following

$$\dot{z}_1(t) = \overline{A}_{11} z_1(t) + \overline{B}_{11} z_1(t - \tau), \tag{8}$$

$$\dot{z}_2(t) = \overline{A}_{22} z_2(t) + \overline{B}_{22} z_2(t - \tau), \tag{9}$$

where $z(t)$ is defined as same as Example 1 and

$$\overline{A}_{11} = \begin{bmatrix} 0 & 1 \\ -1 & 1 \end{bmatrix}, \overline{B}_{11} = \begin{bmatrix} 0 & 0 \\ 0 & 1 \end{bmatrix}, \overline{A}_{22} = \begin{bmatrix} 0 & 0 & -1 \\ 1 & 0 & 1 \\ 1 & -1 & 1 \end{bmatrix}, \overline{B}_{22} = \begin{bmatrix} 0 & 0 & 0 \\ 0 & 0 & 0 \\ 1 & 0 & 0 \end{bmatrix}.$$

By the method [5], system (8) is always unstable while system (9) has two stable poles for $3.1416 < \tau < 3.3077$. For this system the crossing frequencies are $s = \pm j$ and $s = \pm\sqrt{1 + \sqrt{2}}j$. However, the proposed method in [5] confirms that all characteristic roots of the system (7) are in the right half plane, but by decomposing, we understand that the system has two stable poles in a specified time delay interval. The time domain response of the subsystem (8) with constant initial function equals to one, is sketched in Fig. 3. Also, Fig. 4 confirms that the subsystem (9) has two stable poles for $\tau = 3.2$.

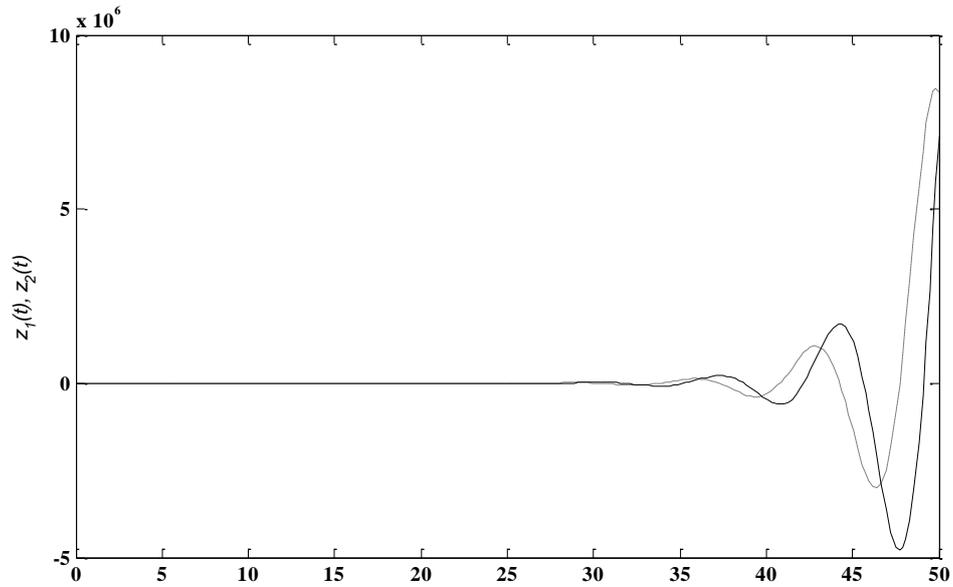

Fig. 3. The time domain response of the system (8). Solid line depicts $z_1(t)$ and dots display $z_2(t)$.

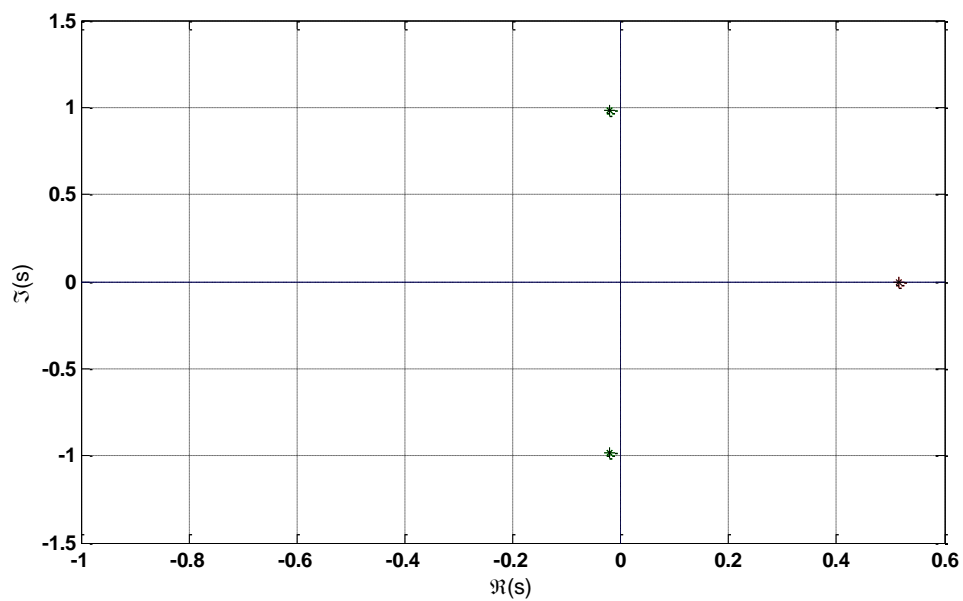

Fig. 4. Roots of the characteristic equation of the system (9), for $\tau = 3.2$.

**4. Stabilization of unstable time delay systems by delayed feedback**

As we said in the introduction, delayed feedback can be employed to stabilize an unstable time-delay system. However, there are infinitely many roots for the characteristic of a retarded DDE, but the number of unstable poles is finite [1]. With considering this issue, we are attempting to move unstable poles of a time-delay system to the left half complex plane by the delayed feedback. Indeed, stabilization is possible, if the characteristic of an unstable time-delay system crosses the imaginary axis. More precisely, for stabilizing by the delayed feedback, there must be exist two crossing frequencies at least, since, the larger crossing frequency corresponds to one where the roots cross from left to right (i.e. destabilizing) of the complex plane [5]. After finding crossing frequencies (if they exist), we can compute the corresponding delays to each crossing frequency. Finally, the interval of delay that the close loop system is stable can be obtained by the cluster treatment [6] or by the direct method [5].

Now, as a concrete example of stabilizing by the delayed feedback, we consider the subsystem (5) that its stability analysis is done in the previous section. Let the open loop time-delay system be as following:

$$\dot{z}(t) = A_1 z(t) + A_2 z(t - 3.2), \tag{10}$$

where $A_1$ and $A_2$ are the same as $\overline{A}_{11}$ and $\overline{B}_{11}$ respectively. It has been shown that the system has two unstable poles.

Now, our goal is the stabilization of the following closed loop system by the delayed feedback:

$$\dot{z}(t) = A_1 z(t) + A_2 z(t - 3.2) - BK\big(z(t) - z(t - \tau)\big), \tag{11}$$

where $B = \begin{bmatrix} 1 \\ 0 \end{bmatrix}$ and $K = [k_1 \quad k_2]$. As we said before, stabilization of (11) requires that the characteristic equation of the system (11) cross the imaginary axis. Before going further, we propose the following lemma that give us a necessary condition such that $s = w_c j$ be a root of the characteristic equation of the system (11).

**Lemma 2.** If $s = \omega_{c1}$, where $\omega_{c1} \in \{\omega \in \mathbb{R} : |w| \leq \beta, \text{for some } \beta > 0\}$, is a root of the characteristic equation of the system (11), then the following relation holds:

$$k_2 - |k_2| - \beta(1 - \beta) \leq 1 - k_1 + |k_1|(3 + \beta)$$

**Proof.** The proof is straightforward. By separating the real and imaginary parts of the characteristic equation of the system (11) at $s = j\omega$, we get the following:

$$1 + \cos(\omega(\tau + 3.2))k_1 - k_1\omega \sin(\omega\tau) + (k_1 + k_2)\cos(\omega\tau) -$$

$$k_1 \cos(3.2\omega) - \omega \sin(3.2\omega) - \omega^2 - k_1 - k_2 = 0$$

$$-\sin(\omega(\tau + 3.22))k_1 - k_1\omega \cos(\omega\tau) - (k_1 + k_2)\sin(\omega\tau) +$$

$$k_1 \sin(3.2\omega) - \omega \cos(3.2\omega) + k_1(\omega - 1) = 0$$

After applying triangle inequality and the famous inequalities $|\sin(x)|, |\cos(x)| \leq 1$, the result follows. □

Now, we come back to the stabilization process. According to lemma 2, if we choose $k_1 = 1$ and $k_2 = -5$, then there are two crossing frequencies and corresponding delays as follows:

$$\omega_{c1} = \pm 1.6564, \tau_{11} = 0.4540, \tau_{12} = 4.2473, \ldots,$$

$$\omega_{c2} = \pm 3.5116, \tau_{21} = 0.9469, \tau_{22} = 2.7362, \ldots.$$

After obtaining the crossing frequencies, there are two scenarios that can be used to get the stability interval. First one is the direct method mentioned in [5]. This method yields at the larger crossing frequency roots moves to the right half plane and the next larger corresponds to stabilizing one. Therefore, the system (11) is stable for $0.4540 < \tau < 0.9469$. In second scenario, we apply the cluster treatment method [6].

The root tendency at the larger crossing frequency is +1 and it is -1 at the next one. So, the time-delay system (11) is stable for $0.4540 < \tau < 0.9469$. In fact, two approach yields the same result. The time domain response and characteristic roots of the system (11) are illustrated in Fig. 5 and Fig. 6 respectively.

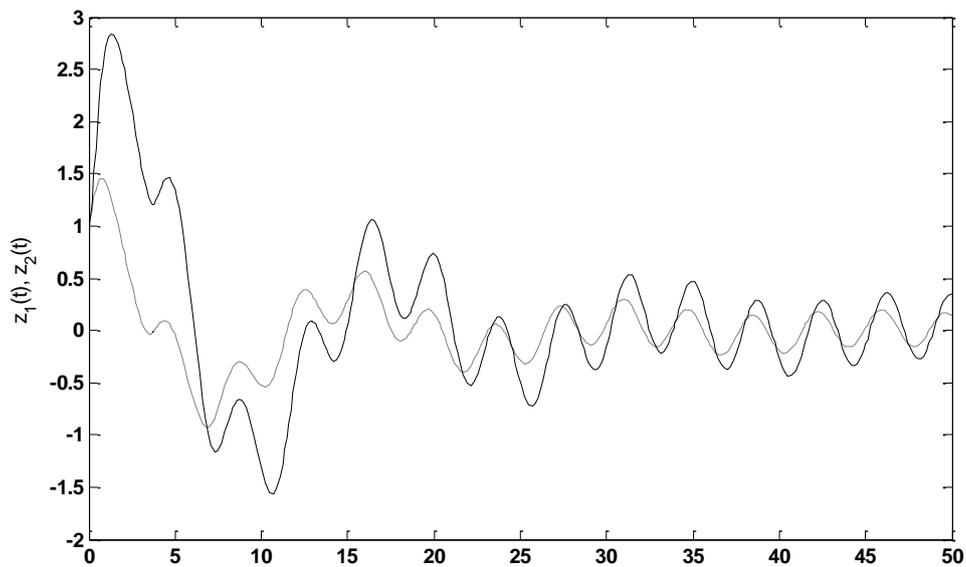

Fig. 5. The time domain response of the system (11) for $\tau = 0.5, k_1 = 1, k_2 = -5$. Solid line depicts $z_1(t)$ and dots display $z_2(t)$.

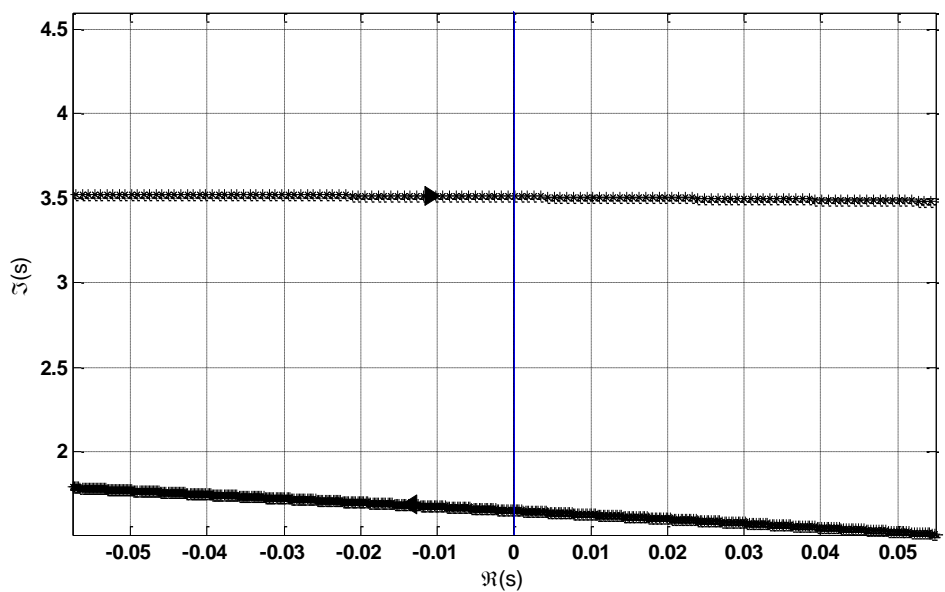

Fig. 6. The root locus of the system (11) near the imaginary axis for $0 < \tau < 2$.

## 5. Design controller via delayed feedback

Now, we apply the delayed feedback to produce a desired response for a DDE. We know that the response of a stable linear time invariant system extremely depends on its dominant poles. It is also true for a DDE, since by the method of steps [1], a DDE convert to an ODE. Summarizing, the response of a stable DDE can be determined by its dominant poles. In general case, i.e. when the DDE has some unstable poles, the rightmost poles determine the response. Here, we consider the subsystem (6) and we change the location of its dominant poles by delayed feedback. Let

$$\dot{z}(t) = A_1 z(t) + A_2 z(t - 3.2), \tag{12}$$

where $A_1$ and $A_2$ are the same as $\overline{A}_{22}$ and $\overline{B}_{22}$ respectively. As we see in Fig. (7), the settling time for this system is very high.

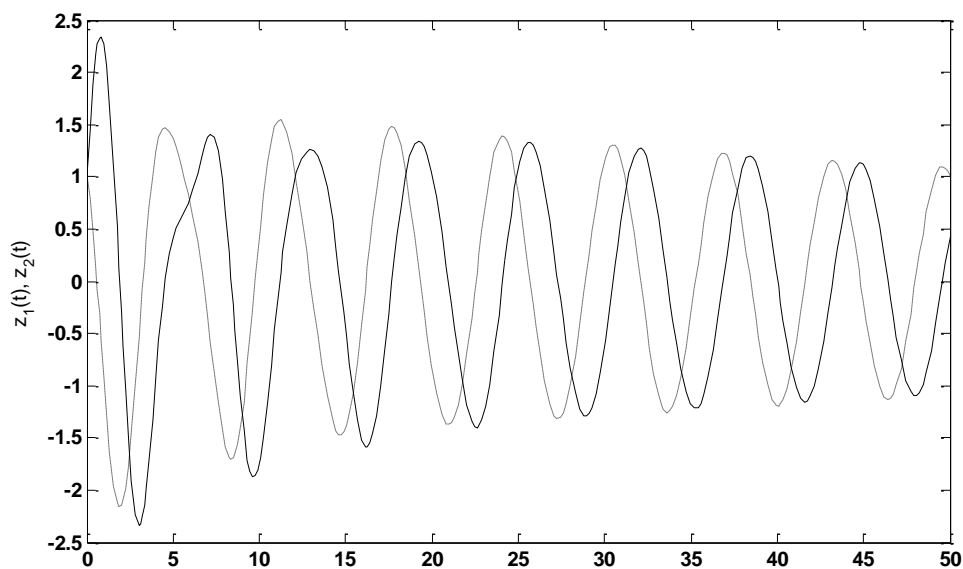

Fig. 7. The time domain response of the system (12). Solid line depicts $z_1(t)$ and dots display $z_2(t)$.

To reduce this time, we consider the following close loop system:

$$\dot{z}(t) = A_1 z(t) + A_2 z(t - 3.2) - BK\big(z(t) - z(t - \tau)\big), \tag{13}$$

where $B = \begin{bmatrix} 1 \\ 0 \end{bmatrix}$ and $K = [k_1 \quad k_2]$. If we force $s = -0.3254 \pm 0.3254j$ to be poles of the system (13), the settling time is highly reduced. After doing this, we obtain two equations. Here, we have three parameters and so we choose $\tau$ freely. If we put $\tau = 0.1$, then we have $k_1 = 40.5925, k_2 = -105.0352$. Finally, the response of the system (13) is plotted in Fig. (8).

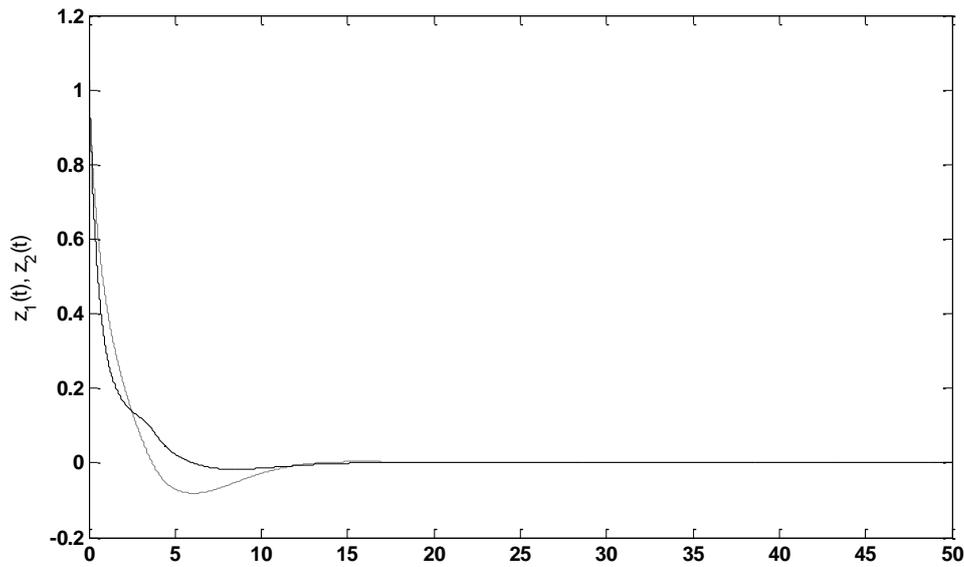

Fig. 8. The time domain response of the system (12) for $\tau = 0.1, k_1 = 40.5925, k_2 = -105.0352$. Solid line depicts $z_1(t)$ and dots display $z_2(t)$.

## 6. Conclusion

In this paper, we established a necessary and sufficient condition for a set of matrices to put them in a block triangular form simultaneously. Also we presented some systems that the stability analysis of them was impossible, unless we decomposed matrices of them to the block form. Furthermore, we employed the delayed feedback to stabilize an unstable time-delay system. More generally, the proposed approach can be applied to an arbitrarily unstable time-delay system. One the other hand, we has been modified the settling time for a time-delay system and so the delayed feedback can also improve the performance.